\documentclass[12pt]{article}

\usepackage{forest}
\usepackage{amsthm, amssymb}
\usepackage{amsmath}               % great math stuff
\usepackage{amsfonts}              % for blackboard bold, etc
\usepackage{amsthm}                % better theorem environments
\usepackage{mathrsfs}

\usepackage{verbatim} %this package enables us to create multiline comments
\usepackage{xspace}
\usepackage{graphicx}  % to include figures
\usepackage[small]{caption} %for fig captions, etc.
\usepackage{color}
\newcommand{\breakingcomma}{%
  \begingroup\lccode`~=`,
  \lowercase{\endgroup\expandafter\def\expandafter~\expandafter{~\penalty0 }}}
\usepackage{float}
\usepackage{wrapfig}

\usepackage{titlesec}
\titleformat{\section}{\large\bfseries}{\thesection}{.5em}{}

\newtheorem{theorem}{Theorem}[section]
\newtheorem{lemma}[theorem]{Lemma}
\newtheorem{remark}[theorem]{Remark}
\newtheorem{proposition}[theorem]{Proposition}
\newtheorem{corollary}[theorem]{Corollary}
\newtheorem{definition}[theorem]{Definition}
\newtheorem{example}[theorem]{Example}
\newcount\refno
\newcommand\be{\begin{equation}}
\newcommand\ee{\end{equation}}
\newcommand\bn{\begin{eqnarray}}
\newcommand\en{\end{eqnarray}}
\newcommand\bns{\begin{eqnarray*}}
\newcommand\ens{\end{eqnarray*}}
\newcommand\bd{\begin{definition}}
\newcommand\ed{\end{definition}}
\newcommand\br{\begin{remark}}
\newcommand\er{\end{remark}}
\newcommand\bt{\begin{theorem}}
\newcommand\et{\end{theorem}}
\newcommand\bp{\begin{proposition}}
\newcommand\ep{\end{proposition}}
\newcommand\bc{\begin{corollary}}
\newcommand\ec{\end{corollary}}
\newcommand\bl{\begin{lemma}}
\newcommand\el{\end{lemma}}

\newcommand\bK{{\mathbb K}}
\newcommand\bR{{\mathbb R}}
\newcommand\bN{{\mathbb N}}

\newcommand\cA{{\cal A}}
\newcommand\cB{{\cal B}}
\newcommand\cD{{\cal D}}
\newcommand\cE{{\cal E}}
\newcommand\cR{{\cal R}}
\newcommand\cL{{\cal L}}
\newcommand\cM{{\cal M}}
\newcommand\cN{{\cal N}}

\newcommand{\N}{\mbox{$\mathbb{N}$}}
\newcommand{\NS}{\mbox{\scriptsize${\mathbb{N}}$}}

\newcommand{\F}{\mbox{$\mathcal{F}$}}

\allowdisplaybreaks

\begin{document}

\title{Sequence Characterization of Multiple Almost-Riordan Arrays and Their Compressions}

\author{Tian-Xiao He 
\\
{\small Department of Mathematics}\\
 {\small Illinois Wesleyan University}\\
 {\small Bloomington, IL 61702-2900, USA}\\
}

\date{}

%%%%%%%%%%%%%%%%%%%%%%%%%%%%%%%%%%%%%%%%%%%%%%%%%%%%%%%%%%%%%%%%%%%%%%
\maketitle
\setcounter{page}{1}
\pagestyle{myheadings}
\markboth{T. X. He }
{Multiple Almost-Riordan Arrays and their sequence characterizations}

%%%%%%%%%%%%%%%%%%%%%%%%%%%%%%%%%%%%%%%%%%%%%%%%%%%%%%%%%%%%%%%%%%%%%%
\begin{abstract}
\noindent 
This is the second paper of the paper series on multiple Riordan arrays. In this paper, based on the study of multiple Riordan arrays and the multiple Riordan group, we define multiple almost-Riordan arrays and find that the set of all multiple almost-Riordan arrays forms a group, called the multiple almost-Riordan group. We also obtain the sequence characteristics of multiple almost-Riordan arrays and give the production matrices for multiple almost-Riordan arrays. We define the compression of multiple almost-Riordan arrays and provide their sequence characterization.



\end{abstract}

\section{Introduction}

Riordan arrays are infinite, lower triangular matrices defined by the generating function of their columns. They form a group, denoted by $\mathcal{R}$ and called {\em the Riordan group} (see Shapiro, Getu, W. J. Woan and L. Woodson \cite{SGWW}). 

A Riordan array $(g(t),f(t))$ is a pair of formal power series $g(t) = \sum_{n\geq 0}g_nt^n$ and $f(t) = \sum_{n\geq 1}f_nt^n$, with $g_0\not= 0$ and $f_1\not= 0$. It defines an infinite lower triangular array $[dn,k]_{n,k\geq 0}$ according to the rule $d_{n,k} = [t^n ]g(t)f(t)^k$. The set of all Riordan arrays forms a group under matrix multiplication

\[
(g(t), f (t))(h(t), l(t)) = (g(t)h(f (t)), l(f (t))).
\]
An almost-Riordan array is defined by an ordered triple $(a|g, f )$ of power series where $a(t) =\sum_{n\geq 0} a_nt^n$, $g(t) = \sum_{n\geq 0} g_nt^n$ and $f(t) = \sum_{n\geq1} f_nt^n$, with $a_0, g_0, f_1\not= 0$. The array is identified with the lower-triangular matrix defined as follows: its first column is given by the expansion of $a(t)$. The remaining element of the infinite tridiagonal matrix coincide with the Riordan array $(g(t), f (t))$.

An infinite lower triangular matrix $D = [d_{n,k}]_{n,k\geq 0}$ is called double Riordan array, if $g(t)$ gives column zero and $f_1(t)$ and $f_2(t)$ are multiplier functions, where $g(t) =\sum_{n\geq 0} g_{2n}t^{2n}$, $f_1(t) =\sum_{n\geq 0} f_{1,2n+1}t^{2n+1}$ and $f_2(t) =\sum_{n\geq 0} f_{2,2n+1}t^{2n+1}$ together with $g_0,f_{1,1}, f_{2,1}\not= 0$. Then the double Riordan array related to $g(t), f_1(t)$ and $f_2(t)$, denoted by $(g; f_1, f_2)$, has the column vector $(g; gf_1, gf_1f_2, gf_1^2f_2, gf_1^2f_2^2, \ldots)$.

Riordan arrays, almost-Riordan arrays and double Riordan arrays have emerged as a powerful tool with broad applications in various branches of mathematics. With their intricate connections to combinatorics, group theory, matrix theory, and number theory, these arrays serve as a bridge between these disciplines. This paper presents the study of the double almost-Riordan arrays and the double almost-Riordan group defined in \cite{He24}. Specifically, it studies the sequence characterization, compression, and total positivity of double almost-Riordan arrays.

More formally, let us consider the set of all formal power series (f.p.s.) in $t$, $\F = {\mathbb K}[\![$$t$$]\!]$, with a field ${\mathbb K}$ of characteristic $0$ (e.g., ${\mathbb Q}$, ${\mathbb R}$, ${\mathbb C}$, etc.). The \emph{order} of $f(t)  \in \F$, $f(t) =\sum_{k=0}^\infty f_kt^k$ ($f_k\in {\bR}$), is the minimal number $r\in\N$ such that $f_r \neq 0$. Denote by $\F_r$ the set of formal power series of order $r$. Let $g(t) \in \F_0$ and $f(t) \in \F_1$. Then, the pair $(g(t) ,\,f(t) )$ defines the {\em (proper) Riordan array} $D=(d_{n,k})_{n,k\in \NS}=(g(t), f(t))$ having
  
\begin{equation}\label{Radef}
d_{n,k} = [t^n]g(t) f(t) ^k
\end{equation}
or, in other words, having $g(t) f(t)^k$ as the generating function whose coefficients make-up the entries of column $k$. 

From the {\it fundamental theorem of Riordan arrays} (see \cite{Sha1}), it is immediate to show that the usual row-by-column product of two Riordan arrays is also a Riordan array:
\begin{equation}\label{Proddef}
    (g_1(t) ,\,f_1(t) )  (g_2(t) ,\,f_2(t) ) = (g_1(t) g_2(f_1(t) ),\,f_2(f_1(t) )).
\end{equation}
The Riordan array $I = (1,\,t)$ acts as an identity for this product. Thus, the set of all Riordan arrays forms the Riordan group $\mathcal{R}$.

Several subgroups of $\mathcal{R}$ are important and have been considered in the literature:
\begin{itemize} \item the set $\mathcal{A}$ of {\em Appell arrays} is the collection of all Riordan arrays $R = (g(t) ,\,t )$ in ${\cR}$; 
\item the set $\mathcal{L}$ of {\em Lagrange arrays} is the collection of all Riordan arrays $R = (1 ,\,f(t) )$ in ${\cR}$;
\item the set $\mathcal{B}$ of \emph{Bell} or {\em renewal arrays} is the collection of all Riordan arrays $R = (g(t) ,\,t g(t))$ in ${\cR}$;
\item the set $\mathcal{H}$ of \emph{hitting-time arrays} is the collection of all Riordan arrays $R = (tf'(t)/f(t) ,\,f(t))$ in ${\cR}$;
\item the set $\mathcal{D}$ of the Riordan arrays $R = (f'(t) ,\, f(t) )$ in ${\cR}$ is called the {\it derivative subgroup}.
\item the set $\mathcal{E}$ of the Riordan arrays $R=((f(t)/t)^{r}f'(t)^{s}, f(t))$ in ${\cR}$ with real or complex $r$ and $s$ is called {\it Lu\'zon-Merlini-Mor\'on-Sprugnoli (LMMS) subgroup}, denoted by ${\cE}[r,s]$, which includes ${\cD}={\cE}[0,1]$ as its special case (see \cite{LMMS2}). 
\end{itemize}

From \cite{Rog}, an infinite lower triangular array $[d_{n,k}]_{n,k\in{\bN}}=(g(t), f(t))$ is a Riordan array if and only if an {\it $A$-sequence} $A=(a_0\not= 0, a_1, a_2,\ldots)$ exists such that for every $n,k\in{\bN}$ there holds 
\be\label{eq:1.1}
d_{n+1,k+1} =a_0 d_{n,k}+a_1d_{n,k+1}+\cdots +a_nd_{n,n},
\ee 
which is shown in \cite{HS} to be equivalent to 
\be\label{eq:1.2}
f(t)=tA(f(t)).
\ee
Here, $A(t)$ is the generating function of the $A$-sequence. In \cite{MRSV} it is also shown that a unique {\it $Z$-sequence} $Z=(z_0, z_1,z_2,\ldots)$ exists such that every element in column $0$ can be expressed as the linear combination 
\be\label{eq:1.3}
d_{n+1,0}=z_0 d_{n,0}+z_1d_{n,1}+\cdots +z_n d_{n,n},
\ee
or equivalently (see \cite{HS}),
\be\label{eq:1.4}
g(t)=\frac{d_{0,0}}{1-tZ(f(t))}.
\ee

Denote the {\it upper shift matrix} by $U$, i.e., 
\[
U=(\delta_{i+1,j})_{i,j\geq 0}=\left [ \begin{array}{llllll} 0& 1& 0& 0& 0& \cdots\\
0&0 &1& 0& 0&\cdots\\
0 &0& 0&1& 0& \cdots\\
0&0& 0& 0&1&\cdots\\
\vdots &\vdots& \vdots& \vdots&\vdots&\ddots\end{array}\right]
\]
and 
\be\label{eq:1.4-2}
P=\left [ \begin{array}{llllll} z_0& a_0& 0& 0& 0& \cdots\\
z_1& a_1 &a_0& 0& 0&\cdots\\
z_2 &a_2& a_1& a_0& 0& \cdots\\
z_3& a_3& a_2& a_1&a_0&\cdots\\
\vdots &\vdots& \vdots& \vdots&\vdots&\ddots\end{array}\right]=\left( Z(t), A(t), tA(t), t^{2}A(t),\ldots\right), 
\ee
where the rightmost expression is the representation of $P$ by using its column generating functions. Here, $P$ is called the {\it production matrix} or {\it $P$-matrix characterization} or simply {\it $P$ matrix} (see Deutsch, Ferrari, and Rinaldi \cite{DFR05, DFR}). From \cite{DFR05, DFR} or Proposition 2.7 of \cite{He15}, the $P$-matrix of Riordan array $R$ satisfies 

\be\label{1.5}
P =R^{-1}UR=R^{-1}\overline{R},
\ee
where $\overline{R}$ is the truncated Riordan array $R$ with the first row omitted. $P$ can be written in terms of $A$- and $Z$-sequences as 

\be\label{eq:1.8}
P=\left [ \begin{array}{llllll} z_0& a_0& 0& 0& 0& \cdots\\
z_1& a_1 &a_0& 0& 0&\cdots\\
z_2 &a_2& a_1& a_0& 0& \cdots\\
z_3& a_3& a_2& a_1&a_0&\cdots\\
\vdots &\vdots& \vdots& \vdots&\vdots&\ddots\end{array}\right],
\ee
where $z_j, a_j\geq 0$ for $j=0,1,\ldots$.

In construction of a Riordan array, one multiplier function h is used to multiply one column to obtain the next column. Suppose alternating rules are applied to generate an infinite matrix similar to a Riordan array. To consider this case, one may use $\ell$, $\ell \geq 2$ multiplier functions, denoted by $f_1$, $f_2\ldots,$ and $f_\ell$, respectively.  

Let $g\in {\mathbb K}[[t^\ell]]$ with $g(t)=\sum_{k\geq 0} g_{\ell k}t^{\ell k}$ and $f\in t{\mathbb K}[[t^\ell]]$ with $f_j(t)=\sum_{k\geq 0} f_{j, \ell k+1}$ $t^{\ell k+1}$, $j=1,2,\ldots, \ell$. Then the multiple Riordan matrix in terms of  $g(t)$, $f_i(t)$, $i=1,\ldots, \ell$, denoted by $(g; f_{1}, f_{2}, \ldots, f_\ell)$, is defined by the generating function of its columns as 

\[
(d_{n,k})_{n,k\geq 0}=(g, gf_{1}, gf_{1}f_{2}, \ldots, gf_1f_2\cdots f_\ell, gf_{1}^{2}f_{2}\cdots f_\ell, 
\ldots),
\]
where 

\[
d_{n,k}=[t^n]g f_1^{\lfloor \frac{k+\ell -1}{\ell}\rfloor}f_2^{\lfloor \frac{k+\ell -2}{\ell}\rfloor}\cdots f_\ell^{\lfloor \frac{k}{\ell}\rfloor}.
\]

Here, we use $\ell$ cases of the {\it first fundamental theorem of multiple Riordan type arrays}:

\be\label{1.6}
(g; f_{1}, f_{2},\ldots, f_\ell)A_j(t)=B_j(t),
\ee
where for $A_j(t)=\sum_{k\geq 0}a_{\ell k+j}t^{\ell k+j}$, $j=0,1,\ldots, \ell-1$,
we have $B_0(t)=gA(h)$ and $B_j(t)=g(f_j/h)A(h)$, $j=1,\ldots, \ell-1$, where $h=\sqrt[\ell]{f_{1}\cdots f_{\ell}}$. Based on the fundamental theorem of multiple Riordan type arrays, we may define a multiplication of two multiple Riordan type arrays as 

\begin{align}\label{1.7}
&(g; f_{1},f_{2}, \ldots f_\ell)(d; h_{1},h_{2},\ldots, h_\ell)\nonumber\\
=&\left(g d(h); \frac{f_1}{h}h_1(h), 
\frac{f_2}{h}h_2(h),\ldots, \frac{f_\ell }{h}h_\ell (h)\right),
\end{align}
where $d(t)$ = $\sum _{k=0}^{\infty} d_{\ell k}t^{\ell k}$, $h_{j}(t)$ = $\sum_{k=0}^{\infty} h_{j,\ell k} t^{\ell k +1}$, $j=1,2,\ldots, \ell$.

Hence, the multiple Riordan array $(g; f_1,\ldots, f_\ell)$ is corresponding to the multiple Riordan type array $(g; f_1/t, \ldots, f_\ell/t)$. 

The collection of all multiple Riordan arrays forms the multiple Riordan group under the multiplication defined by \eqref{1.7}, in which $d(t)$ = $\sum _{k=0}^{\infty} d_{\ell k}t^{\ell k}$, $h_{j}(t)$ = $\sum_{k=0}^{\infty} h_{j,\ell k+1} t^{\ell k+1}$, $j=1,2,\ldots, \ell$. The collection of all multiple Riordan arrays forms a group called the multiple Riordan group and denoted by ${\cal M}\cal {R}$. Clearly, the identity of ${\cal M}{\cal R}$ is $(1;t,\ldots, t)$, and the inverse of $(g; f_1,\ldots, f_\ell)$ is 

\be\label{1.8}
(g; f_1,\ldots, f_\ell)^{-1}=\left(\frac{1}{g(\bar h)}; \frac{t\bar h}{f_1(\bar h)}, \frac{t\bar h}{f_2(\bar h)}, \ldots, \frac{t\bar h}{f_\ell(\bar h)}\right),
\ee
where $\bar h$ is the compositional inverse of $h=\sqrt[\ell]{f_1f_2\cdots f_\ell}.$

Similarly, the bivariate generating function of the multiple Riordan group element $(g; f_1, f_2, \ldots f_\ell)$ is given by 

\be\label{1.6-2-2}
\frac{g(1 + yf_1 + \cdots +y^{\ell-1}f_1f_2\cdots f_{\ell-1})}{1-y^\ell f_1f_2\cdots f_\ell}.
\ee
In particular, the row sums and the diagonal sums of the associated matrix have generating functions

\be\label{1.6-3-2}
\frac{g(1 + f_1 + \cdots +f_1f_2\cdots f_{\ell-1})}{1-f_1f_2\cdots f_\ell}
\ee
and

\be\label{1.6-4}
\frac{g(1 + xf_1 + x2f_1f_2 +\cdots +x^{\ell -1}f_1f_2\cdots f_\ell)}{1-x^\ell f_1f_2\cdots f_\ell},
\ee
respectively.

The study of the case $\ell=2$ for $\ell$-multiple Riordan arrays, called double Riordan arrays, is started from Davenport, Shapiro, and Woodson \cite{DSW12}, followed by the author \cite{He18}, Branch, Davenport, Frankson, Jones, and Thorpe \cite{BDFJT}, Davenport, Frankson, Shapiro, and Woodson \cite{DFSW}, Sun and Sun \cite{SS23}, and Zhang and Zhao \cite{ZZ}, etc. and their references. The cases of $\ell=3$ and $4$ are studied in Barry \cite{Bar24} in a different view.  The set of all double Riordan arrays forms the double Riordan group denoted by ${\cD}{\cR}$. 

\begin{proposition}\label{pro:2.3}
We can identify some subgroups of the multiple Riordan group ${\cal M}{\cal R}$:

\begin{itemize} 
\item the set $\mathcal{A}$ of {\em Appell arrays} is the collection of all multiple Riordan arrays $\{(g; t\ldots, t):g\in \F^{(\ell)}_0\}$ in ${\cal M}{\cR}$, which is a subgroup and called the Appell subgroup of ${\cal M}{\cR}$; 

\item the set $\mathcal{L}$ of {\em Lagrange arrays} is the collection of all multiple Riordan arrays $\{ (1; \,f_1,\ldots, f_\ell):f\in \F^{(\ell)}_1\}$ in ${\cal M}{\cR}$, which is a subgroup and called the Lagrange subgroup of ${\cal M}{\cR}$;

\item the set $\mathcal{D}$ of {\it derivative arrays} is the collection of all multiple Riordan arrays $\{ \left(h'(t), f_1,\ldots, f_\ell \right): f_i\in \F^{(\ell)}_1, i=1,\ldots, \ell\}$ with $h=\sqrt[\ell]{f_1\cdots f_\ell}$  in ${\mathcal M}{\cR}$, which is a subgroup and called the derivative subgroup of ${\mathcal M}{\cal M}{\cR}$.

\item the set $\mathcal{B}_j$ of {\it $j$-th Bell arrays} is the collection of all multiple Riordan arrays $\{ (g; f_1,\ldots, f_\ell): g\in \F^{(\ell)}_0, f_i\in \F^{(k)}_1, i=1,\ldots, \ell\,\, and \,\, f_j=tg\}$ in ${\mathcal M}{\cR}$, which is a subgroup and called the $j$-th Bell subgroup of ${\mathcal M}{\cR}$.
\end{itemize} 
\end{proposition}

In \cite{DSW}, subgroups ${\cA}:=\{(1; f_1,f_2)\in{\cD}{\cR}\}$, ${\cB}_1:=\{(g; tg,f_2)\in{\cD}{\cR}\}$, ${\cB}_2:=\{(g; f_1,tg)\in{\cD}{\cR}\}$, and ${\cN}:=\{ (g;t,t)\in{\cD}{\cR}\}$ are given. Here, ${\cN}$ is a normal subgroup of ${\cD}{\cR}$ and ${\cD}{\cR}$ is the semidirect product of ${\cN}$ and ${\cA}$. 

In \cite{He25}, the following sequence characterization of a multiple Riordan array in ${\cM}{\cR}$ is given. 

\begin{theorem}\label{thm:new-3.0}(\cite{He25})
Let $(d_{n,k})_{n, k\geq 0}=(g;f_1,f_2,\ldots, f_\ell)$ be a multiple Riordan array, and let $A(t)=\sum_{k\geq 0}a_kt^{\ell k}$, $Z_i(t)=\sum_{k\geq 0} z_{i,k} t^{\ell k}$, $i=0, 1,2,\ldots, \ell -1$, be the generating functions of $A$-, $Z_0$-, $Z_1$-, $\ldots$, and $Z_{\ell -1}$-sequences, respectively. Then

\begin{align}\label{3.10}
&A(t)=\frac{t^\ell}{\overline {h}^\ell},\\
&Z_0(t)=  \frac{1}{\overline{h}^\ell}\left( 1-\frac{g_0}{g(\overline{h})}\right)\label{3.11}\\
&Z_m(t)=\frac{1}{\overline{h}^\ell}\left( 1-\frac{g_0f_{1,1}f_{2,1}\cdots f_{m,1}\overline{h}^m}{g(\overline{h})f_1(\overline{h})f_2(\overline{h})\cdots f_m(\overline{h})}\right)\label{3.12},
\end{align}
where $m=1,2,\ldots, \ell-1$, and $\overline{h}$ is the compositional inverse of $h=\sqrt[\ell]{f_1f_2\cdots f_\ell}$.
\end{theorem}

We now present the definition of almost-Riordan arrays. 

\begin{definition}\label{def:1.3}\cite{Barry16}
Let $d, g\in \F_0$ with $d(0), g(0)=1$ and $f\in \F_1$ with $f'(0)=1$. 
We call the following matrix an almost-Riordan array with respect to $d,g,$ and $f$ and denote it by $(d|\,g,f)$:

\be\label{1.13}
(d|\,g,f)=(d, tg, tgf, tgf^2,\cdots),
\ee
where $d$, $tg$, $tgf$, $t^2gf\cdots$, are the generating functions of the $0$th, $1$st, $2$nd, $3$rd, $\cdots$, columns of the matrix $(a|\,g,f)$, respectively. 
\end{definition}

It is clear that $(d|\,g,f)$ can be written as 
\be\label{1.14}
(d|\,g,f)=\left( \begin{matrix} d(0) & 0\\ (d-d(0))/t & (g,f) \end{matrix}\right),
\ee
where $(g,f)=(g,gf, gf^2,\ldots)$, a Riordan array. Particularly, if $d=g$ and $f=t$, then the almost-Riordan array $(d|\,g,f)$ reduces to the Appell-type Riordan array $(g, t)$. 

Barry, Pantelidis, and the author \cite{BHP} present the following operation form for the almost-Riordan group.

\begin{theorem}\label{thm:0.3}\cite{BHP}
The set of all almost-Riordan arrays defined by \eqref{1.13} forms a group, denoted by $a{\cR}$, with respect to the multiplication defined by 
\be\label{1.15}
(a|\,g,f)(b|\,d,h)=\left( \left. a+\frac{tg}{f}(b(f)-1)\right |\,gd(f), h(f)\right),
\ee
where $a, g, b, d\in \F_0$ and $f, h\in \F_1$, and $(a|\,g,f)$ and $(b|\,d,h)$ are the almost-Riordan arrays defined by \eqref{1.13} or \eqref{1.14}.
\end{theorem}

Alp and Kocer give the sequence characterization of almost-Riordan arrays and define the exponential almost-Riordan group in \cite{AK} and \cite{AK24}, respectively. Slowik and the author discuss the total positivity of almost-Riordan arrays in \cite{HS24} and the total positivity of quasi-Riordan arrays in \cite{HS}, respectively, where the set of quasi-Riordan arrays forms a normal subgroup of the almost-Riordan group defined in \cite{He}. 

Let us introduce one more group that is closely related to the groups of Riordan and almost-Riordan arrays.

\begin{definition}\label{def:0.3}\cite{He} 
Let $g\in \F_0$ with $g(0)=1$ and $f\in \F_1$. We call the following matrix a quasi-Riordan array and denote it by $[g,f]$.
\be\label{0.6}
[g,f]:=(g,f,tf,t^2f,\ldots),
\ee
where $g$, $f$, $tf$, $t^2f\cdots$ are the generating functions of the $0$th, $1$st, $2$nd, $3$rd, $\cdots$, columns of the matrix $[g,f]$, respectively. It is clear that $[g,f]$ can be written as 
\be\label{0.8}
[g,f]=\left( \begin{matrix} g(0) & 0\\ (g-g(0))/t & (f/t,t)\end{matrix}\right),
\ee
where $(f,t)=(f,tf,t^2f,t^3f,\ldots)$ and $(f/t, t)$ is an Appell Riordan array. 
Particularly, if $f=tg$, then the quasi-Riordan array $[g,tg]=(g, t)$, a Appell-type Riordan array. 

Clearly, $[g,f]=(g|f/t, t)$. 

Let $A$ and $B$ be $m\times m$ and $n\times n$ matrices, respectively. Then we define the direct sum of $A$ and $B$ by

\be\label{0.5} 
A\oplus B =\left[ \begin{matrix} A &0 \\ 0 &B\end{matrix}\right]_{(m+n)\times(m+n)}.
\ee

In this notation the Riordan array $(g,f)$ satisfies 

\be\label{0.9}
(g,f)=[g,f]([1]\oplus (g,f)).
\ee
\end{definition}

Denote by $q{\cR}$ the set of all quasi-Riordan arrays defined by \eqref{0.6}. In \cite{He} it is shown that $q{\cR}$ is a group with respect to regular matrix multiplication. More precisely, there is the following result.

\begin{theorem}\label{thm:1.2}\cite{He}
The set of all quasi-Riordan arrays $q{\cR}$ is a group, called the quasi-Riordan group, with respect to the multiplication represented in 
\be\label{1.11}
[g,f][d,h]=\left[g+\frac{f}{t}(d-1), \frac{fh}{t}\right],
\ee
which is derived from the first fundamental theorem for quasi-Riordan arrays (FFTQRA),

\be\label{1.12}
[g,f]u=gu(0)+\frac{f}{t}(u-u(0)).
\ee
Hence, $[1,t]$ is the identity of $q{\cR}$. 
\end{theorem}

In \cite{BHP}, it has been shown that the quasi-Riordan group is a normal subgroup of the almost-Riordan group. In Slowik and the author's paper \cite{HS23}, the total positivity of quasi-Riordan arrays is discussed. For instance, let $(g(t), f(t))$ be a Riordan array, where $g(t)=\sum_{n\geq 0} g_nt^n$ and $f(t)=\sum_{n\geq 1}f_nt^n$. If the lower triangular matrix 

\begin{equation}\label{eq:1}
\begin{array}{rl}
Q & =[g,f]=(g|f/t,t)=
\left [ \begin{array}{llllll} g_0& 0& 0& 0& 0& \cdots\\
g_1& f_1 & 0& 0& 0&\cdots\\
g_2 &f_2& f_1& 0& 0& \cdots\\
g_3& f_3& f_2& f_1&0&\cdots\\
\vdots &\vdots& \vdots& \vdots&\vdots&\ddots\end{array}\right]
\\
& =\left( g(t), f(t), tf(t), t^{2}f(t),\ldots\right)
\\
\end{array}
\end{equation}
is totally positive (TP), then so is $R=(g,f)$ (cf. \cite{MMW} or \cite{He, HS23}).  Other interesting criteria for total positivity of Riordan arrays can be found in \cite{CLW,CW}.

In Slowik and the author's another paper \cite{HS24}, a simdirect product is given and used to discuss the total positivity of almost-Riordan arrays via the total positivity of quasi-Riordan arrays.

\begin{theorem}\label{thm:2.5}\cite{HS24}
Every almost-Riordan array $(d|\,g,f)$ can be written as the semidirect product 

\be\label{eq:R_factoriz-1}
(d|\, g,f)=[d,tg](1|\,1,f),
\ee
or equivalently, 

\begin{equation}
\label{eq:R_factoriz}
(d|\,g,f)=\left[\begin{array}{c|c}
d_0 & 0  \\
\hline 
\frac{d-d_0}{t}&  (g,t)\\
\end{array}\right]
\left[\begin{array}{c|c}
1 & 0\\
\hline 
0 & (1,f)\\
\end{array}\right].
\end{equation}
\end{theorem}

The paper will be organized as follows. In next section, we continue the work shown in \cite{He25} to define  multiple almost-Riordan arrays and the multiple almost-Riordan group  from another perspective by using an alternative view. In Section $3$, we give the sequence characterizations of multiple almost-Riordan arrays and the production matrices of multiple almost-Riordan arrays. In Section $4$, we discuss algebraic properties of the multiple almost-Riordan group. Finally, we provide the compression of multiple almost-Riordan arrays and their sequence characterization in Section $5$.

\section{Multiple almost-Riordan arrays}

\begin{definition}\label{def:6.1}
Let $b,g\in {\mathbb K}[[t^\ell]]$ with $b(t)=\sum_{k\geq 0}b_{\ell k}t^{\ell k}$, $g(t)=\sum_{k\geq 0} g_{\ell k}t^{\ell k}$ and $f\in t{\mathbb K}[[t^\ell]]$ with $f_j(t)=\sum_{k\geq 0} f_{j, \ell k+1}$ $t^{\ell k+1}$, $j=1,2,\ldots, \ell$. Then the multiple almost-Riordan matrix in terms of  $b(t)$, $g(t)$, $f_i(t)$, $i=1,\ldots, \ell$, denoted by $(b|g; f_{1}, f_{2}, \ldots, f_\ell)$, is defined by the generating function of its columns as 

\begin{align}\label{6.1}
&(b|g;f_1,f_2,\ldots, f_\ell)=(d_{n,k})_{n,k\geq 0}\nonumber\\
=&(b, tg, tgf_1, tgf_1f_2, \ldots, tgf_1f_2\cdots f_\ell, tg f_1^2f_2\cdots f_\ell, tg f_1^2f_2^2\cdots f_\ell,\ldots),
\end{align}
where 

\begin{align}\label{6.2}
&d_{n,0}=[t^n] b(t)=\left\{ \begin{array}{ll} b_{\ell k} &\mbox{if $n=\ell k$}\\
0 &\mbox{if $n\not= \ell k$}\end{array}\right.\nonumber\\
&d_{n,k}=[t^{n-1}]g f_1^{\lfloor \frac{k+\ell -2}{\ell}\rfloor}f_2^{\lfloor \frac{k+\ell -3}{\ell}\rfloor}\cdots f_\ell^{\lfloor \frac{k-1}{\ell}\rfloor}.
\end{align}
\end{definition}

\begin{theorem}\label{thm:3.2}
There two cases of the first fundamental theorems of multiple almost-Riordan arrays:

\be\label{3.2}
(b|g; f_{1}, f_{2}, \ldots, f_\ell)u(t)=v(t),
\ee
where for $u(t)=\sum_{k\geq 0}u_{\ell k}t^{\ell k}$, $u(t)=\sum_{k\geq 0}u_{2k+1}t^{2k+1}$, and $u(t)=\sum_{k\geq 0}u_{2k+j}t^{2k+j}$, $j=2,\ldots, \ell-1$, we have 

\begin{align}
v(t)=&u_0b+\frac{tg}{f_\ell}(u(h)-u_0)\quad \mbox{and}\label{3.3}\\
v(t)=&\frac{tg}{h}u(h),\label{3.4}\\
v(t)=&\frac{tgf_1f_2\cdots f_{j-1}}{h^j}u(h),\label{3.4-2}
\end{align}
respectively, where $j=2,3,\ldots, \ell$ and $h=\sqrt[\ell]{f_1f_2\cdots f_\ell}$. 
\end{theorem}

By using the first fundamental theorems of multiple almost-Riordan arrays, we may establish a multiplication operator in the set of multiple almost-Riordan arrays. 

\begin{theorem}\label{thm:3.3}
Let $(b|g;f_1,f_2,\ldots, f_\ell)$ and $(c|d;h_1,h_2,\ldots, h_\ell)$ be two multiple  almost-Riordan arrays, where $b$ and $c$ are formal power series in ${\bK}[[t^\ell]]$, $b=\sum_{k\geq 0} b_{\ell k} t^{\ell k}$ and $c=\sum_{k\geq 0} c_{\ell k}t^{\ell k}$, and 
$f_j, h_j\in t{\bK}[[t^\ell]]$, $j=1,2,\ldots, \ell$. Then 

\begin{align}\label{3.6}
&(b|g;f_1,f_2,\ldots, f_\ell)(c|d;h_1,h_2,\ldots, h_\ell)\nonumber\\
=&\left( c_0b+\frac{tg}{f_\ell}(c(h)-c_0)|\,gd(h),  \frac{f_1}{h}h_1(h), \frac{f_2}{h}h_2(h), 
\ldots, \frac{f_\ell}{h}h_\ell (h)\right),
\end{align}
and 

\begin{align}\label{3.7}
&(b|g;f_1,f_2,\ldots, f_\ell)^{-1}\nonumber\\
=&\left( \frac{1}{b_0}+\frac{f_\ell(\overline{h})}{b_0\overline{h}g(\overline{h)}}(b_0-b(\overline{h}))| \frac{1}{g(\overline{h})};
\overline{h}\frac{t}{f_1(\overline{h})},\overline{h}\frac{t}{f_2(\overline{h})}, \ldots, \overline{h}\frac{t}{f_\ell (\overline{h})}\right),
\end{align}
where $b_0\not= 0$ because $b\in \F_0$, and $\overline{h}$ is the compositional inverse of $h=\sqrt[\ell]{f_1f_2\cdots f_\ell}$, i.e., 

\[
\overline{h}(\sqrt[\ell]{f_1f_2\cdots f_\ell})=\sqrt[\ell]{f_1f_2\cdots f_\ell}(\overline{h})=t.
\] 
\end{theorem}

\section{Sequence characterization of multiple almost-Riordan arrays}

We now give a sequence characterization of a multiple almost-Riordan array in ${\cM}a{\cR}$. Inspired by Branch, Davenport, Frankson, Jones, and Thorpe \cite{BDFJT} and  Davenport, Frankson, Shapiro, and Woodson \cite{DFSW}, we consider $D=(b|g;f_1,f_2)$ as 

\begin{align}\label{4.9}
D=&(b, tg, tgf_1, tg(f_1f_2), tgf_1(f_1f_2), tgf_1f_2(f_1f_2\cdots f_\ell, \ldots)\nonumber\\
=&(b,0,0,\ldots)+(0,tg, 0, \ldots, tg(f_1f_2\cdots f_\ell), 0, \ldots, tg(f_1f_2\cdots f_\ell)^2, 0, \ldots)\nonumber\\
&+(0,0,tgf_1, 0, \ldots, tgf_1(f_1f_2\cdots f_\ell), 0, \ldots, tgf_1(f_1f_2\cdots f_\ell)^2, 0\ldots)
\nonumber\\
=&D_0+D_1+D_2+\cdots +D_\ell.
\end{align}
After omitting zero columns and top zero rows, we denote the remaining $D_j$, $j=1,2,\ldots, \ell$, shown above by $D_j^*$, $j=1,2,\ldots, \ell$, respectively. Then $D_1^*=(g, f_1f_2\cdots f_\ell)=(d^{(1)}_{n,k})_{n,k\geq 0}$, $D_2^*=(gf_1/t, f_1f_2\cdots f_\ell)=(d^{(2)}_{n,k})_{n,k\geq 0}, \ldots$, and $D_j^*=(gf_1f_2\cdots f_{\ell-1}, f_1f_2\cdots f_\ell)=(d^{(\ell)})_{n,k\geq 0}$ are Riordan arrays. Hence, $(b|g;f_1,f_2, \ldots, f_\ell)$ has a $W$-sequence, $\ell$ $Z$-sequences, denoted by $Z_j$-sequences, $j=1,2,\ldots, \ell$, and a $A$-sequence in this view, while using the view shown in \cite{He24}, we have $W$-sequence, a $Z$-sequnce, and $\ell$ $A$-sequences: $A_j$-sequences, $j=1,2,\ldots, \ell$. 

By using \eqref{4.9}, we have the following result.

\begin{theorem}\label{thm:new-3.1}
Let $(b|g;f_1,f_2, \ldots, f_\ell)$ be a multiple almost-Riordan array, and let $A(t)=\sum_{k\geq 0}a_kt^{\ell k}$, $Z_j(t)=\sum_{k\geq 0} z_{j,k} t^{\ell k}$, $j=1,2,\ldots, \ell$, and $W(t)=\sum_{k\geq 0} w_k t^{\ell k}$ be the generating functions of $A$-, $Z_j$-, $j=1,2,\ldots, \ell$, and $W$-sequences, respectively. Then

\begin{align}\label{8.0}
&A(t)=\frac{t^\ell }{{\overline h}^\ell},\\
&Z_m(t)=\left\{ \begin{array}{ll} \frac{1}{\overline{h}^\ell}\left( 1-\frac{g_0}{g(\overline{h})}\right) & \mbox{if $m=1$}\\
\frac{1}{\overline{h}^\ell}\left( 1-\frac{g_0f_{1,1}f_{2,1}\cdots f_{m-1,1}\overline{h}^{m-1}}{g(\overline{h})f_1(\overline{h})f_2(\overline{h})\cdots f_{m-1}(\overline{h})}\right) &\mbox{if $m=2,3,\ldots, \ell-1$}\end{array}, \right. \label{8.1}\\
&Z_\ell(t)=  \frac{g_0f_{1,1}f_{2,1}\cdots f_{\ell-1,1}}{b_0}+\frac{t^\ell}{\overline{h}^\ell}-\frac{g_0f_{1,1}f_{2,1}\cdots f_{\ell -1,1}b(\overline{h})f_\ell(\overline{h})}{b_0\overline{h}g(\overline{h})},\label{8.2}\\
&W(t)=\frac{f_\ell(\overline{h})((1-w_0\overline{h}^\ell)b(\overline{h})-b_0)}{\overline{h}^{\ell+1}g(\overline{h})}+w_0,\quad w_0=b_\ell/b_0,\label{8.3}
\end{align}
where $m=2,3,\ldots, \ell-1$, $f_{j,1}=[t] f_j(t)$, $j=1,2,\ldots, \ell$, and $\overline{h}$ is the compositional inverse of $h=\sqrt[\ell]{f_1f_2\cdots f_\ell}$.
\end{theorem}

We may using production matrix to represent the sequences characterizations of $(b|g;f_1,f_2)$ in terms of the view shown in \cite{BDFJT, DFSW,He25}. 

\begin{theorem}\label{thm:4.3-2}
Let $(b|g;f_1,f_2,\ldots, f_\ell)\in {\cM}a{\cR}$. Then, $(b|g;f_1,f_2,\ldots, f_\ell)$ has a production matrix 

\be\label{4.2-2}
P=\left( W(t), tZ_1(t), \ldots, t^{\ell-1}Z_{\ell -1}(t), Z_\ell(t), tA(t), t^2A(t), t^3A(t), \ldots \right),
\ee
where $A(t)$, $Z_j(t)$, $j=1,2,\ldots, \ell$, and $W(t)$ are shown in \eqref{8.0}-\eqref{8.3}.
\end{theorem}

An alternative sequence characterization and the corresponding production matrix for multiple almost-Riordan arrays can be found by using the view shown in Theorem $3.1$ of \cite{He25}. 

\begin{example}\label{ex:2.1}
Consider Example $1$ given in Barry \cite{Bar24} for the case $\ell=3$, we may set 
a triple almost-Riordan array as follows.

\begin{align*}
&b(t)=\frac{1}{1-t^6},\\
&g(t)=\frac{1}{1-t^3},\\
&f_1(t)=\frac{t}{1-t^3},\\
&f_2(t)=t(1+t^3),\\
&f_3(t)=\frac{t}{1+t^3}.
\end{align*}
Thus 

\[
h(t)=\sqrt[3]{f_1f_2f_3}=\frac{t}{\sqrt[3]{1-t^3}}\quad \mbox{and}\quad \overline {h}(t)=\frac{t}{\sqrt[3]{1+t^3}}.
\]
The first few rows of Riordan type array $(g;f_1,f_2,f_3)$ is 

\[
\left(\frac{1}{1-t^6}|\frac{1}{1-t^3};\frac{t}{1-t^3},t(1+t^3),\frac{t}{1+t^3}\right)=\left [ \begin{array}{lllllllllll} 
1&0& 0& 0& 0& 0&0& 0& 0& 0&  \cdots\\
0&1& 0& 0& 0& 0&0& 0& 0& 0&  \cdots\\
0&0&1&0 &0& 0& 0& 0& 0& 0& \cdots\\
0&0&0& 1&0& 0& 0& 0& 0& 0& \cdots\\
0&1&0& 0& 1&0& 0& 0& 0& 0& \cdots\\
0&0&2& 0& 0&1& 0& 0& 0& 0& \cdots\\
1&0&0& 3& 0&0& 1& 0& 0& 0& \cdots\\
0&1&0& 0& 2&0& 0& 1& 0& 0& \cdots\\
0&0&3& 0& 0&3& 0& 0& 1& 0& \cdots\\
0&0&0& 5& 0&0& 4& 0& 0& 1& \cdots\\
\vdots &\vdots& \vdots& \vdots&\vdots&\vdots &\vdots &\vdots &\vdots &\ddots\end{array}\right].
\]

From Theorem \ref{thm:new-3.1}, 

\begin{align*}
&A(t)=1+t^3,\\
&Z_1(t)=1,\\
&Z_2(t)=1+\frac{1}{1+t^3}=1-t^3+t^6-t^9+\cdots,\\
&Z_3(t)=2+t^3-\frac{(1+t^3)^2}{(1+2t^3)^2}=1+3t^3-5t^6+12t^9+\cdots,\\
&W(t)=\frac{t^3(1+t^3)}{(1+2t^3)^2}=t^3-3t^6+8t^9+\cdots.
\end{align*}
Hence, $A$-sequence, $Z_1$-, $Z_2$-, $Z_3$-, and $W$-sequences are 

\begin{align*}
&A=(1,1,0,\ldots), \\
&Z_1=(1,0,\ldots),\\
&Z_2=(2,-1, 1,-1,\ldots),\\
&Z_3=(1,3,-5,12,\ldots),\\
&W=(0,1, -3,8,\ldots).
\end{align*} 

For instances, 

\begin{align*}
&d_{6,0}=w_0d_{3,0}+w_1d_{3,3}=1,\\
&d_{4,1}=z_{1,0}d_{1,1}=1,\\
&d_{5,2}=z_{2,0}d_{2,2}=1,\,\, d_{8,2}=z_{2,0}d_{5,2}+z_{2,1}d_{5,5}=3,\\
&d_{3,3}=z_{3,0}d_{0,0}=1,\,\,d_{6,3}=z_{3,0}d_{3,0}+z_{3,1}d_{3,3}=3,\\
&d_{9,3}=z_{3,0}d_{6,0}+z_{3,1}d_{6,3}+z_{3,2}d_{6,6}=5,\\
\end{align*}

The production matrix of $(g;f_1,f_2, f_3)$ is 

\[
P=\left [ \begin{array}{lllllllllll} 
0&0& 0& 0& 0& 0&0& 0& 0& 0&  \cdots\\
0&1& 0& 0& 1& 0&0& 0& 0& 0&  \cdots\\
0&0& 2& 0& 0& 1& 0& 0& 0& 0& \cdots\\
1&0& 0& 1& 0& 0& 1& 0& 0& 0& \cdots\\
0&0& 0& 0& 1&0& 0& 1& 0& 0& \cdots\\
0&0&-1& 0& 0&1& 0& 0& 1& 0& \cdots\\
-3&0& 0& 3& 0&0& 1& 0& 0& 1& \cdots\\
0&0& 0& 0& 0&0& 0& 1& 0& 0& \cdots\\
0&0& 1& 0& 0&0& 0& 0& 1& 0& \cdots\\
8&0& 0& -5& 0&0& 0& 0& 0& 1& \cdots\\
\vdots &\vdots& \vdots& \vdots&\vdots&\vdots &\vdots &\vdots &\vdots &\ddots\end{array}\right].
\]
\end{example}

\begin{example}\label{ex:2.2} 
Considering the double almost-Riordan array, $(1/(1-t^4)|1/(1-t^2); t, t/(1-t^2))$, in which $(1/(1-t^2);t, t/(1-t^2))$ presents the Fibonacci-Stanley tree. Here, $f_1=t$, $f_2=t/(1-t^2)$, $g=1/(1-t^2)$, and $b=1/(1-t^4)$.  Thus, $h=\sqrt{f_1f_2}=t/\sqrt{1-t^2}$, and the compositional inverse of $h$ is $\overline{h}=t/\sqrt{1+t^2}$. Substituting $f_1=t$ $f_2=t/(1-t^2)$, $g=1/(1-t^2)$, and $b=1/(1-t^4)$ into Equations \eqref{8.0}-\eqref{8.3} for $\ell =2$ and noting 

\begin{align*}
&f_1(\overline{h})=\frac{t}{\sqrt{1+t^2}},\\
&f_2(\overline{h})=t\sqrt{1+t^2},\\
&g(\overline{h})=1+t^2,\\
&b(\overline{h})=\frac{(1+t^2)^2}{1+2t^2},
\end{align*}
we immediately have 

\begin{align*}
&A(t)=\frac{t^2}{\overline{h}^2}=1+t^2,\\
&Z_1(t)=\frac{1}{(t/\sqrt{1+t^2})^2}\left(1-\frac{1}{1+t^2}\right)=1,\\
&Z_2(t)=\frac{t\sqrt{1+t^2}\left( 1+t^2-\frac{(1+t^2)^2}{1+2t^2}\right)}{\frac{t}{\sqrt{1+t^2}}(1+t^2)}+1=\frac{1+3t^2+t^4}{1+2t^2},\\
&\qquad = 1+t^2-t^4+2t^6-4t^8+8t^{10}+\cdots,\\
&W(t)=\frac{t\sqrt{1+t^2}\left(\frac{(1+t^2)^2}{1+2t^2}-1\right)}{\left(\frac{t}{\sqrt{1+t^2}}\right)^3(1+t^2)}+0=\frac{t^2(1+t^2)}{1+2t^2}\\
&\qquad =t^2-t^4+2t^6-4t^8+8t^{10}+\cdots.
\end{align*}
The double almost-Riordan array $(1/(1-t^4)|1/(1-t^2); t, t/(1-t^2))$ begins 

\be
\left[\begin{array}{lllllllllll}
1&0&0&0&0&0&0&0&0&0&...\\
0&1&0&0&0&0&0&0&0&0&...\\
0&0&1&0&0&0&0&0&0&0&...\\
0&1&0&1&0&0&0&0&0&0&...\\
1&0&1&0&1&0&0&0&0&0&...\\
0&1&0&2&0&1&0&0&0&0&...\\
0&0&1&0&2&0&1&0&0&0&...\\
0&1&0&3&0&3&0&1&0&0&..\\
1&0&1&0&3&0&3&0&1&0&...\\
0&1&0&4&0&6&0&4&0&1&...\\
\vdots&\vdots&\vdots&\vdots&\vdots&\vdots&\vdots&\vdots&\vdots&\vdots&\ddots
\end{array}\right]
\ee

Denote 

\begin{align*}
P=&(W(t), tZ_1(t), Z_2(t), tA(t), t^2A(t), \ldots)\\
=&\left[\begin{array}{llllllllllll}
0& 0& 1&0&0&0&0&0&0&0&0&...\\
0& 1& 0&1&0&0&0&0&0&0&0&...\\
1& 0& 1&0&1&0&0&0&0&0&0&...\\
0& 0 &0&1&0&1&0&0&0&0&...\\
-1&0&-1&0&1&0&1&0&0&0&0&...\\
0& 0&0&0&0&1&0&1&0&0&0&...\\
2& 0& 2&0&0&0&1&0&1&0&0&...\\
0&0&  0&0&0&0&0&1&0&1&0&..\\
-4&0&-4&0&0&0&0&0&1&0&1&...\\
0& 0& 0&0&0&0&0&0&0&1&0&...\\
\vdots&\vdots&\vdots&\vdots&\vdots&\vdots&\vdots&\vdots&\vdots&\vdots&\ddots
\end{array}\right].
\end{align*}

We find 

\begin{align*}
(b|g;f_1.f_2)P=
&\overline{\overline{\left( \frac{1}{1-t^4}|\frac{1}{1-t^2}; t,\frac{t}{1-t^2}\right)}},
\end{align*}
where the rightmost matrix is the truncation of 

\[
(b|g;f_1,f_2)=\left( \frac{1}{1-t^4}|\frac{1}{1-t^2}; t, \frac{t}{1-t^2}\right)
\]
with the first row and the second element of column zero omitted. 
\end{example}

\begin{example}\label{ex:2.3} Kuznetkov, Pak,  and Postnikov \cite{KPP} consider ordered trees with no points at odd heights, which is converted to the Dyck paths with no valley at odd heights, but ending at the heights in \cite{DSW}. Denote by $d_{n,k}$ the number of those paths with $n$ edges that end at height $k$. Then the matrix $D=(d_{n,k})_{n,k\geq 0}$ is 

\be\label{2.5-3}
D=
\left[\begin{array}{lllllllllll}
1&0&0&0&0&0&0&0&0&0&...\\
0&1&0&0&0&0&0&0&0&0&...\\
0&0&1&0&0&0&0&0&0&0&...\\
0&1&0&1&0&0&0&0&0&0&...\\
0&0&2&0&1&0&0&0&0&0&...\\
0&2&0&2&0&1&0&0&0&0&...\\
0&0&4&0&3&0&1&0&0&0&...\\
0&4&0&5&0&3&0&1&0&0&..\\
0&0&9&0&8&0&4&0&1&0&...\\
0&9&0&12&0&9&0&4&0&1&...\\
\vdots&\vdots&\vdots&\vdots&\vdots&\vdots&\vdots&\vdots&\vdots&\vdots&\ddots\\
\end{array}\right]
\ee
It can be seen that 

\begin{align*}
&d_{n,2k}=d_{n-1,2k-1}+d_{n-1,2k+1},\quad k\geq 1,\\
&d_{n, 2k+1}=d_{n-2, 2k-1}+d_{n-2,2k+1}+d_{n-2,2k+3}, \quad k\geq 1,\\
&d_{n,1}=d_{n-1, 0}+d_{n-1,2},\\
&d_{n,0}=0,
\end{align*}
for $n\geq 1$. Thus, $A_1=1+t^2$, $A_2=1+t^2+t^4$, $Z=1+t^2$, and $W=0$.
\end{example}

\section{Algebraic properties of the multiple almost-Riordan group}
We now discuss some subgroups of the multiple almost-Riordan arrays. A multiple almost-Riordan array $(b|g;f_1,f_2,\ldots, f_\ell)$ is said to be normalized, if $b(0)=g(0)=1$. 

\begin{theorem}\label{thm:5.1}
Let ${\cD}=(b|g; f_1,f_2,\ldots, f_\ell)\in {\cM}a{\cR}$ be normalized. If $f_j=t$, $j=1,2,\ldots, \ell$, then the collection of all elements $(b|g;t,t,\ldots, t)$ in ${\cM}a{\cR}$, denoted by ${\cA}$, forms a subgroup, called the Appel subgroup of ${\cM}a{\cR}$, i.e., ${\cM}\leq {\cM}a{\cR}$.  

If $b=g=1$, then the collection of all elements $(1|1;f_1,f_2,\ldots, f_\ell)$ in ${\cM}a{\cR}$, denoted by ${\cL}$, forms a subgroup, called the Lagrange subgroup of ${\cM}a{\cR}$, i.e., ${\cL}\leq {\cM}a{\cR}$. 

For $j=1,2,\ldots, \ell$, if $f_j=tg$, then the collection of all elements $(b|g;tg,$ $f_2,\ldots, f_\ell)$ in ${\cM}a{\cR}$, denoted by ${\cB_j}$, forms a subgroup, called the type-j Bell subgroup of ${\cM}a{\cR}$, i.e., ${\cB_j}\leq {\cD}a{\cR}$.
\end{theorem}

\begin{theorem}\label{thm:5.2}
The Appell subgroup ${\cA}$ is a normal subgroup of the multiple almost-Riordan group, i.e., ${\cA} \vartriangleleft {\cM}a{\cR}$. There exists the following semidirect  product for ${\cM}a{\cR}$: 

\be\label{3.8}
{\cM}a{\cR}={\cA}\rtimes{\cL}
\ee
\end{theorem}

The following theorem shows how to extend some multiple Riordan subgroups to multiple almost-Riordan subgroups. 

\begin{theorem}\label{thm:5.4}
Let $(b|g;f_1,f_2,\ldots, f_\ell)\in{\cM} a{\cR}$, and let $S$ be a subgroup of ${\cM}{\cR}$. Then $\{ (b|g;f_1,f_2,\ldots, f_\ell): b\in \F_0, (g;f_1,f_2,\ldots, f_\ell)\in S\}$ be a subgroup of ${\cM}a{\cR}$. Therefore, ${\mathcal D}=\{ (b|h';f_1,f_2, \ldots, f_\ell):h=\sqrt[\ell]{f_1f_2\cdots f_\ell}\}$ is a subgroup of ${\cM}a{\cR}$. ${\mathcal B_j}=\{ (b|g;f_1,f_2,\ldots, f_\ell): f_j=tg\}$, $j=1,2,\ldots, \ell$, is a subgroup of ${\cM}a{\cR}$.
\end{theorem}

The following theorem gives a new subgroup of ${\cal M}a{\cal R}$ and shows how to extend some Riordan subgroups to multiple almost-Riordan subgroups by using the new subgroup of ${\cM}a{\cR}$. 

\begin{theorem}\label{thm:5.5}
For $j=1,2,\ldots, \ell$, $\widehat {\cL_j}:=\{(b|g;f_1,f_2,\ldots, f_\ell)\in{\cM}a{\cR}, f_j=t\}$ is a subgroup of ${\cM}a{\cR}$. Furthermore, for $\ell=2$, if $\{(g,f_2)\in {\cR}\}$ is in Appel subgroup or type-$2$ Bell subgroup of ${\cR}$, then corresponding $\{ (b|g;t,f_2)\in{\cD}a{\cR}\}$ is in Appel subgroup or type-$2$ Bell subgroup of ${\cD}a{\cR}$, respectively. 
\end{theorem}

\section{Compressions of multiple almost-Riordan arrays} 

Let $(b|g;f_1,f_2,\ldots, f_\ell)=(d_{n,k})_{n\geq k\geq 0}\in {\cal M}a{\cal R}$. We define its compression $(\hat d_{n,k})_{n\geq k\geq 0}$ as follows:

\be\label{9.1} 
\hat d_{n,k}:=d_{\ell n-(\ell -1)k,k}, \quad n\geq k\geq 0. 
\ee
We now study the structure of the compression of a Riordan array starting from the following theorem. 

\begin{theorem}\label{thm:9.2}
Let $(b|g;f_1,f_2,\ldots, f_\ell)=(d_{n,k})_{n\geq k\geq 0}$ be a multiple almost-Riordan array with 

\begin{align*}
& b(t)=\sum_{k\geq 0} b_{k}t^{\ell k}, \quad g(t)=\sum_{k\geq 0} g_{k}t^{\ell k},\nonumber\\ 
&f_j(t)=\sum_{k\geq 0}f_{j,k}t^{\ell k+1}, \,\,  j=1,2,\ldots, \ell,
\end{align*}
and let its compression array $(\hat d_{n,k})_{n,k\geq 0}$ be defined by \eqref{9.1}. Then we have  

\begin{align}\label{9.1-2} 
&\hat d_{n,0}=[t^n] \hat b(t),\nonumber \\
&\hat d_{n,k}=\begin{cases} [t^n]t\hat g (\hat f_1\hat f_2\cdots \hat f_\ell)^{(k-1)/\ell}, & \mbox{if $k \equiv 1\, (mod\, \ell)$},\\
[t^n]t\hat g \hat f_1\cdots \hat f_{m-1}(\hat f_1\hat f_2\cdots \hat f_\ell)^{(k-m)/\ell}, &\mbox{if $k\equiv m\,(mod\, \ell)$}
\end{cases}
\end{align}
for $k\geq 1$, where 
\begin{align}\label{9.1-3}
&\hat b(t)=\sum_{k\geq 0} b_{2k}t^k, \quad \hat g(t)=\sum_{k\geq 0} g_{2k}t^k,\nonumber\\
&\hat f_j(t)=\sum_{k\geq 0}f_{j,2k+1}t^{k+1}, \,\, j=1,2,\ldots, \ell.
\end{align}
\end{theorem}

\begin{example}\label{ex:4.1}
As an example, $(1/(1-t^2)|1/(1-t), t/(1-t), t(1+t), \frac{t}{1+t})$ is the compression of the triple almost-Riordan array shown in Example \ref{ex:2.1}, while $(1/(1-t^2)|1/(1-t), t, t/(1-t))$ is the compression of the double almost-Riordan array shown in Example \ref{ex:2.2}. 

The first few rows of the compression of $\left(\frac{1}{1-t^6}|\frac{1}{1-t^3};\frac{t}{1-t^3},t(1+t^3),\frac{t}{1+t^3}\right)$ shown in Example \ref{ex:2.1} is 

\begin{align}\label{0.0-1}
&\left(\frac{1}{1-t^2}|\frac{1}{1-t};\frac{t}{1-t},t(1+t),\frac{t}{1+t}\right)\nonumber\\
=&\left [ \begin{array}{llllllllll} 
1&0& 0& 0& 0& 0&0& 0&0&  \cdots\\
0&1& 0& 0& 0& 0&0& 0&0&  \cdots\\
1&1&1&0 &0& 0& 0& 0&0& \cdots\\
0&1&2& 1&0& 0& 0& 0& 0& \cdots\\
1&1&3& 3& 1&0& 0& 0& 0& \cdots\\
0&1&4& 5& 2&1& 0& 0& 0& \cdots\\
1&1&5& 7& 3&3& 1& 0& 0& \cdots\\
0&1&6& 9& 4&6& 4& 1& 0&\cdots\\
1&1&7& 11& 5&10& 9& 4& 1&  \cdots\\
\vdots &\vdots& \vdots& \vdots&\vdots&\vdots &\vdots  &\ddots\end{array}\right].
\end{align}
\end{example}

The sequence characterization of the compression of a multiple almost-Riordan array is given in the following theorem. 

\begin{theorem}\label{thm:9.3}
Let $(b|g;f_1,f_2)=(d_{n,k})_{n\geq k\geq 0}$ be a multiple almost-Riordan array, and let its compression be defined by $(\hat d_{n,k})_{n\geq k\geq 0}$, where $\hat d_{n,k}$ is shown in \eqref{9.1}. 
Suppose the $A$-, $Z_j$-, $j=1,2,\ldots, \ell$, and $W$-sequences of $(b|g;f_1,f_2)$ are 

\begin{align*}
&A=\{ a_0, a_1,\ldots\}, \,\,Z_j=\{ z_{j,0}, z_{j,1}, \ldots\}, \quad j=1,2,\ldots, \ell,
 \mbox{and}\\
 & W=\{w_0,w_1,\ldots\}
\end{align*}
with their generating functions $A(t)=\sum_{k\geq 0}a_kt^{\ell k}$, $Z_j(t)=\sum_{k\geq 0} z_{j,k} t^{\ell k}$, $j=1,2,\ldots, \ell$, and $W(t)=\sum_{k\geq 0} w_k t^{\ell k}$ in ${\bK}[[t^\ell]]$, respectively. Then 

\begin{align}\label{9.2-2}
&A\left( \sqrt[\ell]{\frac{\hat f_1\hat f_2\cdots \hat f_\ell}{t^{\ell -1}}}\right)=\frac{\hat f_1\hat f_2\cdots \hat f_\ell}{t^\ell},\\
&Z_1\left(\sqrt[\ell]{\frac{\hat f_1\hat f_2\cdots \hat f_\ell}{t^{\ell -1}}}\right)=\frac{1}{t}\left(1-\frac{g_0}{\hat g}\right), \label{9.3-2}\\
&Z_m\left(\sqrt[\ell]{\frac{\hat f_1\hat f_2\cdots \hat f_\ell}{t^{\ell -1}}}\right)=\frac{1}{t}\left(1-\frac{g_0f_{1,1}f_{2,1}\cdots f_{m-1,1}t^{m-1}}{\hat g\hat f_1\hat f_2\cdots \hat f_m}\right), \quad m=2,3,\ldots, \ell-1, \label{9.4-2}\\
&Z_\ell\left(\sqrt[\ell]{\frac{\hat f_1\hat f_2\cdots \hat f_\ell}{t^{\ell -1}}}\right)=\frac{\hat f_\ell}{t^\ell \hat g}(\hat g\hat f_1\hat f_2\cdots \hat f_{\ell -1}-z_{\ell,0}t^{\ell -1}\hat b)+z_{\ell,0},\nonumber\\
&\qquad z_{\ell,0}=\frac{g_0f_{11}f_{2,1}\cdots f_{\ell-1, 1}}{b_0},\label{9.4-2-2}\\
&W\left(\sqrt{\frac{\hat f_1\hat f_2\cdots \hat f_\ell}{t^{\ell -1}}}\right)=\frac{\hat f_\ell}{t^2\hat g}\left(\hat b(1-w_0t)-b_0\right)+w_0. \quad w_0=\frac{b_\ell}{b_0},\label{9.5-2}
\end{align}
where $f_{j,1}=[t]f_j$, $j=1,2,\ldots, m$, or equivalently,

\begin{align}\label{9.2} 
\hat d_{n,k}=&a_{0}\hat d_{n-\ell,k-\ell}+a_{1}\hat d_{n-1,k}+a_{2}\hat d_{n+\ell-2,k+\ell}+\cdots\nonumber\\
=&\sum_{j\geq 0} 
a_j\hat d_{n-\ell+j(\ell -1), k+\ell(j-1)},\quad k> \ell\\
\hat d_{n,\ell}=&z_{\ell, 0}\hat d_{n-\ell,0}+z_{\ell,1}\hat d_{n-1, \ell}+z_{\ell,2}\hat d_{n+\ell-2, 2\ell}+\cdots,
\nonumber\\
=&\sum_{j\geq 0}z_{\ell, j}\hat d_{n-\ell+j(\ell -1), j\ell}, \label{9.2-0}\\
\hat d_{n,m}=&z_{m,0}\hat d_{n-1,m}+z_{m,1}\hat d_{n+\ell-2, \ell+m}+z_{m,2} \hat d_{n+2\ell-3, 2\ell +m}+\cdots, \nonumber\\
=&\sum_{j\geq 0}z_{m,j}\hat d_{n-1+j(\ell -1), j\ell+m},\quad m=2,3,\ldots, \ell-1, \label{9.3}\\
\hat d_{n,1}=&z_{1,0}\hat d_{n-1,1}+z_{1,1}\hat d_{n+\ell-2, \ell+1}+z_{1,2}\hat d_{n+2\ell -3, 2\ell +1}+\cdots \nonumber\\
=&\sum_{j\geq 0} z_{1,j}\hat d_{n+j(\ell -1)-1, j\ell +1}, \label{9.4}\\
\hat d_{n,0}=&w_0\hat d_{n-1, 0}+w_1\hat d_{n+\ell -2, \ell}+w_2\hat d_{n+2\ell -3, 2\ell}+\cdots, \nonumber\\
=&\sum_{j\geq 0} w_j\hat d_{n+j(\ell-1)-1, j\ell}, \label{9.5},
\end{align}
\end{theorem}

\begin{remark}\label{rem:4.1}
It can be seen that the compositional inverse of $\sqrt{(\hat f_1\hat f_2\cdots \hat f_\ell)/t^{\ell-1}}$ is

\[
\overline {\sqrt{\frac{\hat f_1 \hat f_2\cdots \hat f_\ell}{t^{\ell-1}}}}=\overline{h}^\ell,
\]
where $\overline{h}$ is the compositional inverse of $h= \sqrt{f_1f_2\cdots f_\ell}$. Hence, \eqref{9.2-2}-\eqref{9.5-2} are equivalent to \eqref{8.0}-\eqref{8.3} correspondingly. For instance, substituting $t=\overline{h}^\ell$ into the expressions involving $A$, $Z_1$, and $Z_m$, $m=2,3,\ldots, \ell-1$, and noting $\hat g(\overline{h}^\ell)=g(\overline{h})$ and 

\be\label{0.0-0}
\hat f_j(\overline{h}^\ell)=\overline{h}^{\ell -1}f_j(\overline{h}), \quad j=1,2,\ldots,\ell,
\ee
we obtain 

\begin{align*}
&A(t)=\frac{t^\ell}{\overline{h}^\ell},\\
&Z_1(t)=\frac{1}{\overline{h}^\ell}\left( 1-\frac{g_0}{g(\overline{h})}\right),\\
&Z_m(t)=\frac{1}{\overline {h}^\ell }
\left( 1-\frac{g_0f_{1,1}f_{2,1}\cdots f_{m-1,1}\overline{h}^{\ell (m-1)}}
{g(\overline{h})f_1(\overline{h})f_2(\overline{h})\cdots f_\ell(\overline{h})\overline{h}^{(m-1)(\ell -1)}}\right)\\
=&\frac{1}{\overline {h}^\ell }\left( 1-\frac{g_0f_{1,1}f_{2,1}\cdots f_{m-1,1}\overline{h}^{m-1}}{g(\overline{h}f_1(\overline{h})f_2(\overline{h})\cdots f_\ell(\overline{h})\overline{h}^{(m-1)(\ell -1)}}\right),
\end{align*}
respectively. To transform $Z_\ell$ from \eqref{9.4-2-2} to \eqref{8.2}, we use $t=\overline{h}^\ell$ and \eqref{0.0-0} 
into \eqref{9.3-2} and notice $f_1(\overline{h})f_2(\overline{h})\cdots f_\ell(\overline{h})=t^\ell$ to obtain

\begin{align*}
&Z_\ell\left(t\right)=\frac{\hat f_1(\overline{h}^\ell)\hat f_2(\overline{h}^\ell)\cdots \hat f_\ell (\overline{h}^\ell)}{\overline{h}^{\ell^2}}-z_{\ell,0}\hat b(\overline{h}^\ell)\frac{\hat f_\ell(\overline{h}^\ell)}{\overline{h}^{\ell}\hat g(\overline{h}^\ell)}+z_{\ell,0}\\
=&\frac{t^\ell}{\overline{h}^\ell}-\frac{z_{\ell,0}f_\ell(\overline{h})b(\overline{h})}{\overline{h} g(\overline{h})}+z_{\ell,0},
\end{align*}
where 

\[
z_{\ell,0}=\frac{g_0f_{11}f_{2,1}\cdots f_{\ell-1, 1}}{b_0}.
\]
Hence, we get \eqref{8.2} from \eqref{9.4-2-2}.

Finally, by applying \eqref{0.0-0} and substituting $t=\overline{h}^2$ into \eqref{9.5-2}, we have 

\begin{align*}
W\left(t\right)=&\frac{\hat f_\ell (\overline{h}^\ell)}{\overline{h}^{2\ell}\hat g(\overline{h}^\ell)}\left(\hat b(\overline{h}^\ell)(1-w_0\overline{h}^\ell)-b_0\right)+w_0\\
=&\frac{f_\ell (\overline{h})}{\overline{h}^{\ell +1}g(\overline{h})}\left(b(\overline{h})(1-w_0\overline{h}^\ell)-b_0\right)+w_0,
\end{align*}
where $w_0=\frac{b_\ell}{b_0}$. Hence, we get \eqref{8.3} from \eqref{9.5-2}. 
\end{remark}

%%%Equation  \eqref{9.2} can be presented 

%%%\bns
%%%&&[t^{n}]t\hat g \hat f_{1}(\hat f_{1}\hat f_{2})^{k -1}=\sum_{i\geq 0}a_{1,i}[t^{n-i-1}]t\hat g(\hat f_{1}\hat f_{2})^{k-1+i}\\
%%%&=&[t^{n}]t\hat g(\hat f_{1}\hat f_{2})^{k -1}\sum_{i\geq 0}a_{1,i}t^{i+1}(\hat f_{1}\hat f_{2})^{i}.
%%%\ens
%%%The above expression has the following equivalent form with the using of the generating function of the columns of the compression matrix $(\hat d_{n,k})_{n,k\geq 0}$:

%%%\[
%%%\hat f_{1}=tA\left( t \hat f_{1}\hat f_{2}\right).
%%%\]

\begin{example}\label{ex:4.2}
$(1/(1-t^2)|1/(1-t), t/(1-t), t(1+t), \frac{t}{1+t})$ is the compression of the triple almost-Riordan array shown in Example \ref{ex:2.1} with its first few rows shown in \eqref{0.0-1}.

From the $A$-sequence, $Z_j$- ($j=1,2,3$), and $W-$-sequences for the triple almost-Riordan array $(1/(1-t^6)|1/(1-t^3), t/(1-t^3), t(1+t^3), \frac{t}{1+t^3})$ shown in Example \ref{ex:2.1}, we obtain the $A$-sequence, $Z_j$- ($j=1,2,3$), and $W$-sequences for the compression triple almost-Riordan array $(1/(1-t^2)|1/(1-t), t/(1-t), t(1+t), \frac{t}{1+t})$:

\begin{align*}
&A=(1,1,0, \ldots), \\
&Z_1=(1,0,\ldots),\\
&Z_2=(2,-1, 1, -1,\ldots),\\
&Z_3=(1,3,-5,12,\ldots),\\
&W=(0,1, -3,8,\ldots). 
\end{align*} 
Hence, we have 

\begin{align*}
&\hat d_{n,k}=\hat d_{n-3,k-3}+\hat d_{n-1,k},\quad k> 3,\\
&\hat d_{n,1}=\hat d_{n-1,1},\\
&\hat d_{n,2}=2\hat d_{n-1,2}-\hat d_{n+1,5}+\hat d_{n+3,8}-\hat d_{n+5, 11}+\cdots,\\
&\hat d_{n,3}=\hat d_{n-3, 0}+3\hat d_{n-1,3}-5\hat d_{n+1,6}+12\hat d_{n+3,9}+\cdots,\\
&\hat d_{n,0}=\hat d_{n+1, 3}-3\hat d_{n+3, 6}+8\hat d_{n+5,9}+\cdots.
\end{align*}

$\hat R$ shown below is the compression of the double Riordan array, the Fibonacci-Stanley array, $(1/(1-t^4)|1/(1-t^2); t, t/(1-t^2))$, given in Example \ref{ex:2.2}, in which $(1/(1-t^2);t, t/(1-t^2))$ presents the Fibonacci-Stanley tree. Namely, 

\[
\hat R=(\hat d_{n,k})_{n,k\geq 0}= 
\begin{bmatrix}
1&0&0&0&0&0&0&0&0&...\\
0&1&0&0&0&0&0&0&0&...\\
1&1&1&0&0&0&0&0&0&...\\
0&1&1&1&0&0&0&0&0&...\\
1&1&1&2&1&0&0&0&0&...\\
0&1&1&3&2&1&0&0&0&...\\
1&1&1&4&3&3&1&0&0&...\\
 \vdots&\vdots&\vdots&\vdots&\vdots&\vdots&\vdots&\vdots&\vdots&\ddots\\
\end{bmatrix} 
\]
Since the $A$-, $Z_1$-, $Z_2$-, and $W$- sequences of $(1/(1-t^4)|1/(1-t^2); t, t/(1-t^2))$ are 
\begin{align*}
&A(t)=\frac{t^2}{\overline{h}^2}=1+t^2,\\
&Z_1(t)=\frac{1}{(t/\sqrt{1+t^2})^2}\left(1-\frac{1}{1+t^2}\right)=1,\\
&Z_2(t)=\frac{t\sqrt{1+t^2}\left( 1+t^2-\frac{(1+t^2)^2}{1+2t^2}\right)}{\frac{t}{\sqrt{1+t^2}}(1+t^2)}+1=\frac{1+3t^2+t^4}{1+2t^2},\\
&\qquad = 1+t^2-t^4+2t^6-4t^8+8t^{10}+\cdots,\\
&W(t)=\frac{t\sqrt{1+t^2}\left(\frac{(1+t^2)^2}{1+2t^2}-1\right)}{\left(\frac{t}{\sqrt{1+t^2}}\right)^3(1+t^2)}+0=\frac{t^2(1+t^2)}{1+2t^2},\\
&\qquad =t^2-t^4+2t^6-4t^8+8t^{10}+\cdots,
\end{align*}
we have  $\hat d_{n,k}=\hat d_{n-2, k-2}+\hat d_{n-1,k}$, $k>2$, $\hat d_{n,1}=\hat d_{n-1,1}$, and 

\begin{align*}
&\hat d_{n,2}=\hat d_{n-2,0}+\hat d_{n-1,2}-\hat d_{n,4}+\cdots\\
&\hat d_{n,0}=\hat d_{n,2}-\hat d_{n+1,4}+\cdots.
\end{align*}
\end{example}

\medbreak

\end{document}